\documentclass{amsart}

\usepackage{eufrak}
\usepackage{amssymb,amsmath, amscd, mathrsfs, yfonts, bbm}
\title[Monotonic Independence over a C$^*$-algebra]
{A Combinatorial Approach to  Monotonic Independence over a
C$^*$-algebra}

\author{Mihai Popa}
\address{${}^1$Indiana University at Bloomington,
 Department of Mathematics, Rawles Hall,
 931 E 3rd St, Bloomington, IN 47405}
\email{mipopa@indiana.edu}
\address{${}^2$Institute of Mathematics, Romanian Academy,
 P.O.Box 1-764, Bucharest, RO-70700, Romania}
\DeclareMathAlphabet{\mathpzc}{OT1}{pzc}{m}{it}
\newtheorem{claim}{}[section]
\newtheorem{defn}[claim]{Definition}
\newtheorem{thm}[claim]{Theorem}
\newtheorem{lemma}[claim]{Lemma}
\newtheorem{remark}[claim]{Remark}
\newtheorem{prop}[claim]{Proposition}
\newtheorem{cor}[claim]{Corollary}

\newcommand{\I}{I}
\newcommand{\orIn}{\lfloor \I,n \rfloor}

\newcommand{\gB}{\mathfrak{B}}

\newcommand{\cE}{\mathcal{E}}
\newcommand{\gA}{\mathfrak{A}}
\newcommand{\mE}{\mathcal{E}^m}
\newcommand{\cL}{\mathcal{L}}
\newcommand{\g}{\gamma}
\newcommand{\gh}{\mathfrak{h}}

\newcommand{\lra}{\longrightarrow}

\newcommand{\wme}{\cE^{wm}}
\newcommand{\mul}{Mul[[\gB]]}

\newcommand{\adm}{\mathfrak{a}}
\newcommand{\m}{\frac{m}{2}}
\newcommand{\inv}{\langle -1\rangle}
\newcommand{\gH}{\mathfrak{H}}
\newcommand{\gK}{\mathfrak{K}}

\usepackage{amssymb}
\input xy
\xyoption{all}

\begin{document}

 \maketitle
\bibliographystyle{alpha}
\begin{abstract}
The notion of monotonic independence is considered in a more general
frame, similar to the construction of operator-valued free
probability. The paper presents constructions for maps with similar
properties to the $H$ and $K$ transforms from the literature, semi
inner-product bimodule analogues for the monotone and weakly
monotone product of Hilbert spaces, an ad-hoc version of the Central
Limit Theorem, an operator-valued arcsine distribution as well as a
connection to operator-valued conditional freeness.
\end{abstract}
\section{Introduction}
 An important notion in non-commutative probability is monotonic
independence, introduced by P. Y. Lu and Naofumi Muraki. Since its
beginning, the study of this notion of independence  was done by
constructions, techniques and developments similar to the theory of
free probability. In a 1998 paper published in Mem. AMS (ref.
\cite{speicher1}), R. Speicher developed an operator-valued analogue
of free independence. The present paper addresses similar problems
to ones discussed in \cite{speicher1}, but in the context of
monotonic independence.

Other motivation is that while for the Free Fock space over a
Hilbert space there is a straightforward analogous semi-inner
product bimodule construction, (as illustrated in \cite{pimsner} and
\cite{speicher1})), there are no similar constructions for various
of its deformations, such as the $q$-Fock spaces (\cite{EP}). As
shown in Section \ref{fock}, the monotone and weakly monotone
Fock-like spaces, which are strongly connected to monotonic
independence, admit analogous semi-inner product bimodules.

 The paper is structured in six sections. The
second section presents the definition of the monotonic independence
over an algebra. In the third section there are constructed maps
with similar properties to the maps $H$ and $K$ from the theory of
monotonic independence, as introduced in \cite{muraki1} and
\cite{bercovici1}.The forth section deals with semi-inner product
bi-module analogues of the monotone and weakly monotone products of
Hilbert spaces and algebras of annihilation operators, as introduced
in \cite{muraki1}, \cite{muraki2}, \cite{janusz}. The fifth section
presents a Central Limit Theorem in the frame of monotonic
independence over a C$^*$-algebra and a positivity result concerning
it.  Since in the scalar-valued case the density of the limit
distribution is the arcsine function (\cite{lu2}, \cite{muraki1}),
the limit in Theorem 5.3 can be regarded as an "operator-valued
arcsine law". The last section introduces a notion of conditionally
free product of conditional expectations extending the definition
and positivity results from \cite{mlotk1} and shows a the connection
to monotonic independence analogous to Proposition 3.1 from
\cite{Uwe1}.
%
%
%
%
%
\section{Preliminaries}
Let $\gB$ be an algebra (not necessarily unital). Within this paper,
the notation $\gB_+\langle\xi_1,\dots,\xi_n\rangle$
 will stand for the free
noncommutative algebra generated by $\gB$ and the symbols
$\xi_1,\dots, \xi_n$. For the smaller algebra
$\gB_+\langle\xi_1,\dots,\xi_n\rangle\ominus\gB$ we will use the
notation $\gB\langle\xi_1,\dots,\xi_n\rangle$.

If $\gB$ is a $\ast$-algebra, we can consider $\ast$-algebra
structures on $\gB_+\langle\xi\rangle$ and $\gB\langle\xi\rangle$
either by letting $(\xi)^*=\xi$ (i.e. the symbol $\xi$ is
self-adjoint) or considering $\gB_+\langle\xi, \xi^*\rangle$ with
$(\xi)^*=\xi$.

 We need to also consider an extended notion of
 non-unital complex algebra.
 $\gA$ will be
called a $\gB$-algebra
 if  $\gA$ is an algebra such that $\gB$ is a subalgebra of $\gA$
 or there is an algebra $\widetilde{\gA}$ containing $\gB$ as a
 subalgebra such that $\widetilde{\gA}=\gA\sqcup\gB$. (The symbol
 $\sqcup$ stands for disjoint union).

A map
  $\Phi:\gA\lra\gB$ is said to be  \emph{$\gB$-linear}
 if
\[ \Phi(b_1xb_2 +y)=b_1\Phi(x)b_2 + \Phi(y)
\]
\noindent for all $ x,y\in\gA$  and $ b_1,b_2\in\gB$.

If $\gB$ is a subalgebra of $\gA$ and $\Phi(b)=b$ for all $b\in\gB$,
then $\Phi$ will be called a \emph{conditional expectation}.

\begin{defn}\label{defmonindep}\emph{Suppose that $\gA$ is a
$\gB$-algebra and  $I$ is a totally ordered set.}

 \emph{A family $\{\gA_j\}_{j\in \I}$ of subalgebras of $\gA$
 is said to be monotonically independent over $\gB$ if given
  $X_j\in\gA_j$($ j\in I$), the following conditions are
 satisfied:}
 \begin{enumerate}
\item[(a)]\emph{for all
 $i<j>k$ in $I$ and $A,B\in\gA$:}
   $\Phi(A X_i X_j X_k B)=\Phi(A X_i \Phi(X_j) X_k B)$
\item[(b)]\emph{for all $i_m>\dots>i_1<k_1<\dots<k_n$ in $I$:}
\begin{eqnarray*}
\Phi(X_{i_m}\cdots X_{i_1} &=&
\Phi(X_{i_m})\cdots\Phi( X_{i_1})\\
\Phi(X_{k_1}\cdots X_{k_n})
 &=&
 \Phi(X_{k_1})\dots \Phi(X_{k_n})
 \\
\Phi(X_{i_m}\cdots X_{i_1}X_{k_1}\cdots X_{k_n})
 &=&
 \Phi(X_{i_m})\dots\Phi( X_{i_1})\Phi(X_{k_1})\dots \Phi(X_{k_n})
\end{eqnarray*}
\end{enumerate}
 \emph{The elements $\{X_j\}_{j\in \I}$ from $\gA$ are said to be monotonically independent over
   $\gB$ if the subalgebras of $\gA$ generated by $X_j$ and $\gB$, are
monotonically independent over $\gB$.}

\end{defn}

 Following \cite{muraki1}, or \cite{muraki2}, one may consider the
 stricter definition of monotonic independence replacing the first
 condition with
 \begin{enumerate}
 \item[(a$^\prime$)]$X_iX_jX_k=X_i\Phi(X_j)X_k$ whenever $i<j>k$.
 \end{enumerate}
 \noindent Yet, definition \ref{defmonindep} (similar to \cite{Uwe1})
 suffices for the results within this paper.

\section{The maps $\kappa$, $\rho$ and $\mathfrak{h}$}

Two important instruments in monotonic probability are the maps
$H_X$ and $K_X$ associated to a selfadjoint element $X$ from a
unital $\ast$-algebra $\mathcal{A}$ with a $\mathbb{C}$-linear
functional $\varphi$ such that $\varphi(1)=1$. Namely $H_X$ is
reciprocal Cauchy transform $H_X(z)=(G_X(z))^{-1}$, where $G_X$ is
the Cauchy transform corresponding to $X$:

\[G_X(z)=\varphi((z-X)^{-1})\]
\noindent and the map $K_X$ is given by

\begin{equation*}
K_X(z)=\frac{\psi_X(z)}{1+\psi_X(z)},\ \text{where}\
\psi_X(z)=\varphi\left(zX(1-zX)^{-1}\right).
\end{equation*}

\noindent Their key properties (see \cite{bercovici2}, \cite{Uwe1})
are that for $X,Y$, respectively $U-1, V$ monotonically independent
with respect to $\varphi$, one has:
\begin{eqnarray*}
H_{X+Y}&=&H_X\circ H_Y\\
K_{UV}&=&K_{VU}=K_U\circ K_V.
\end{eqnarray*}

In the scalar-valued case, the moment generating series of $X$, can
be recovered from $H$ and $K$. For the $\gB$-valued setting, the
$n$-th moment of $X$ is the multilinear function
$m_{X,n}:\gB^{n-1}\lra\gB$,
\[
 m_{X,n}(b_1,\dots,b_{n-1})=\Phi(Xb_1X\cdots Xb_{n-1}X).
\]
The mathematical object replacing the moment generating series is a
multilinear function series over $\gB$ (see \cite{dykema}), that
cannot be recovered from a $\gB$-valued analytic map.

The use of analytic tools (such as the Cauchy transform) is strongly
impaired by the previous considerations, hence the combinatorial
approach is very convenient is the present framework. We will first
construct the $\gB$-valued analytic functions $\mathfrak{h}$
replacing $H$, $\kappa$ and $\rho$ replacing $K$. Based on these
constructions, the second part of the section will address the more
general framework of multiplicative function series over an algebra.

In this section we require $\gA$ to be a $\ast$-algebra, $\gB$ to be
a $C^*$-algebra with norm $||\cdot ||$, and $\Phi:\gA\lra\gB$ to be
a positive conditional expectation. $\overline{\gA}$ will denote the
closure of $\gA$ in the topology given by $X\mapsto||\Phi(X^*X)||$.
For simplicity, we will denote the continuous extension of $\Phi$ to
$\overline{\gA}$ also with $\Phi$.
\begin{defn}
\emph{For $X\in\gA$, consider the $\gB$-valued function $\gh_X$}
\[
\{z\in\gB: ||z||<||X||^{-1}\}\ni z
\mapsto
\gh_X(z)=\left(1-z\Phi(X)\right)^{-1}z\in\gB.
\]
\emph{Observe that $\gh$ is an analytic function defined in a
neighborhood of $0\in\gB$ and $\gh(0)=0$.}
\end{defn}

\begin{thm}\label{thm1}If $X,Y\in\gA$ are monotonically independent, then:
\begin{eqnarray*}
\gh_{X+Y}(z)&=&\gh_X\circ \gh_Y(z)\\
\end{eqnarray*}
for $z$ in a neighborhood of $0\in\gB$.
\end{thm}
\begin{proof}First, note that, for $X_1, X_2 \in\gA$
of sufficiently small norm, we have

\begin{equation}\label{eq1}
\sum_{n=0}^{\infty}(X_1+X_2)^n=
\sum_{p=0}^{\infty}\left(\left(\sum_{k=0}^\infty
X_2^k\right)X_1\right)^p\left(\sum_{m=0}^{\infty}(X_2)^m\right).
\end{equation}

Indeed
\begin{eqnarray*}
\sum_{n=0}^{\infty}(X_1+X_2)^n
&=&
 \sum_{m=0}^{\infty}
 \sum_{\alpha_0,\beta_m\geq 0}
\sum_{\alpha_j,\beta_j\geq 1
}X_1^{\alpha_0}X_2^{\beta_0}\dots
X_1^{\alpha_m}X_2^{\beta_m}\\
&=&
 \sum_{n=0}^{\infty}\left(\sum_{\beta_j\geq 0}
\left(\prod_{j=0}^{n}X_2^{\beta_j}X_1\right)\right)
\left(\sum_{m=0}^{\infty}(X_2)^m\right)\\
&=&
 \sum_{p=0}^{\infty}\left(\left(\sum_{k=0}^\infty
X_2^k\right)X_1\right)^p \left(\sum_{m=0}^{\infty}(X_2)^m\right).
\end{eqnarray*}

 Substituting $X_1=zX$ and $X_2=zY$, (\ref{eq1}) becomes:

 \begin{eqnarray*}
\sum_{n=0}^{\infty}(z(X+Y))^n
&=&
 \sum_{p=0}^{\infty}\left(\left(\sum_{k=0}^\infty
(zY)^k\right)zX\right)^p \left(\sum_{m=0}^{\infty}(zY)^m\right),
 \end{eqnarray*}
 and therefore
 \begin{equation*}
(1-z(X+Y))^{-1}z
=
  \sum_{p=0}^{\infty}\left(\left(\sum_{k=0}^\infty
(zY)^k\right)zX\right)^p \left(\sum_{m=0}^{\infty}(zY)^m\right)z.
\end{equation*}
We deduce that
\begin{equation*}
\Phi\left((1-z(X+Y))^{-1}z\right) = \Phi\left(
\sum_{p=0}^{\infty}\left(\left(\sum_{k=0}^\infty
(zY)^k\right)zX\right)^p \left(\sum_{m=0}^{\infty}(zY)^m\right)z
 \right).
 \end{equation*}
 hence
\[
 \gh_{X+Y}(z)
 =
 \Phi
 \left(
\sum_{p=0}^{\infty}\left[\left((1-zY)^{-1}z\right)X\right]^p
 (1-zY)^{-1}z
\right)
\]

Let $Z=(1-zY)^{-1}z\in\overline{\gA}$. $Z$ is in the closure of the
algebra generated by $Y$ and $\gB$. If $X,Y$ are monotonically
independent over $\gB$, the continuity of $\Phi$
  and Definition \ref{defmonindep}(a),
 imply

 \begin{eqnarray*}
 \Phi\left((ZX)^pZ\right)
 &=&
 \Phi(ZXZ\cdots ZXZ)\\
 &=&
\Phi\left(ZX\Phi(Z)\cdots \Phi(Z)XZ\right)
 \end{eqnarray*}
 Since $X\Phi(Z)X\cdots \Phi(Z)X$ is in the algebra generated by $X$
 and $\gB$, \ref{defmonindep}(b) gives
 \[
  \Phi\left((ZX)^pZ\right)
  =
  \left[\Phi(Z)X\right]^p\Phi(Z)
  \]
Therefore
 \begin{eqnarray*}
 \gh_{X+Y}(z)
&=&
 \Phi \left( \sum_{p=0}^{\infty}\left(\gh_Y(z)X\right)^p \gh_Y(z)
\right)\\
&=& \left(\gh_X\circ \gh_Y\right)(z),
 \end{eqnarray*}
as claimed.
\end{proof}

\begin{defn}\emph{For $X\in\gA$ and $z$ in a neighborhood of $0\in\gB$, define the
maps:}
\begin{eqnarray*}
\vartheta_X(z) &=&
\Phi\left((1-zX)^{-1}zX\right)\\
\kappa_X(z)
&=&(1+\vartheta_X(z))^{-1}\vartheta_X(z)\\
\varrho_X(z) &=&\Phi(Xz(1-Xz)^{-1})\\
\rho_X(z)&=&\varrho_X(z)(1+\varrho_X(z))^{-1}.\\
\end{eqnarray*}
\emph{Observe that the above maps are  $\gB$-valued analytic maps
for which $0$ is a fixed point.}
\end{defn}

\begin{thm}\label{thm2}
Let $U, V\in\gA$ be such that $U-1$ and $V$ are monotonically
independent over $\gB$. Then, for $z$ in some neighborhood of
$0\in\gB$,
\begin{eqnarray*}
\kappa_{VU}(z)&=&\left(\kappa_U\circ\kappa_{V}\right)(z)\\
\rho_{UV}(z)&=&\left(\rho_{U}\circ\rho_{V}\right)(z).
\end{eqnarray*}

\end{thm}

\begin{proof} With the notation $U-1=X$, we obtain
\begin{eqnarray*}
\vartheta_{VU}(z)&=& \Phi\left((1-zVU)^{-1}zVU\right)\\
&=&
\Phi((1-zVU)^{-1}zVU)\\
&=&
 \Phi\left(\sum_{k=0}^{\infty}(zVU)^{k}zVU\right)\\
&=& \Phi\left(\sum_{k=0}^{\infty}[zV(X+1)]^{k}zVU\right)\\
&=& \Phi\Big(\sum_{k=0}^{\infty} \sum_{ \substack{
\alpha_1+\dots\alpha_p=k+1\\
\alpha_j\geq 1 }}
(zV)^{\alpha_1}X(zV)^{\alpha_2}X\dots (zV)^{\alpha_p}U\Big).\\
\end{eqnarray*}

\noindent As in the proof of \ref{thm1}, using Definition
\ref{defmonindep}, the above equation becomes

\begin{eqnarray*}
\vartheta_{VU}(z) &=& \Phi \left(\sum_{k=0}^{\infty}
(\vartheta_V(z)X)^k\vartheta_V(z)U\right)\ \\
&=& \Phi\left( (1-\vartheta_V(z)X)^{-1}\vartheta_V(z)U\right)
\\
&=& \Phi\left((1+\vartheta_V(z)-\vartheta_V(z)U)^{-1}
(1+\vartheta_V(z))(1+\vartheta_V(z))^{-1}\vartheta_V(z)U\right)\\
&=& \Phi\left(
\left[(1+\vartheta_V(z))^{-1}(1+\vartheta_V(z)-\vartheta_V(z)U)\right]^{-1}
(1+\vartheta_V(z))^{-1}\vartheta_V(z)U\right)\\
&=&
\Phi\left([1-(1+\vartheta_V(z))^{-1}\vartheta_V(z)U]^{-1}\kappa_V(z)U\right)\\
&=&
\Phi\left([1-\kappa_V(z)U]^{-1}\kappa_V(z)U\right)\\
&=& \vartheta_U(\kappa_V(z))
\end{eqnarray*}
Therefore:
\begin{eqnarray*}
\kappa_{VU}(z)
 &=&
 [1+\vartheta_{VU}(z)]^{-1}\vartheta_{VU}(z)\\
 &=&
[1+\vartheta_{U}(\kappa_V(z))]^{-1}\vartheta_{U}(\kappa_V(z))\\
&=& \kappa_U(\kappa_V(z)).
\end{eqnarray*}

The identity for $\rho$ follows analogously.
\end{proof}

The proofs of Theorems \ref{thm1} and \ref{thm2} do not use the
analyticity of the maps $\gh$, $\kappa$, $\rho$, but only properties
from Definition \ref{defmonindep} and some combinatorial identities
that are true for any formal series. This leads to a easy
reformulation of the results in the more general frame, presented in
\cite{dykema}, of multilinear function series over an algebra.

In the following paragraphes we will briefly remaind the reader the
construction and several results on multilinear function series.

Let $\gB$ be an algebra. We set $\widetilde{\gB}$ equal to $\gB$ if
$\gB$ is unital and to the unitalization of $\gB$ otherwise. For
$n\geq1$, we denote by $\mathcal{L}_n(\gB)$ the set of all
multilinear mappings
\[\omega_n:\underbrace{\gB\times\dots\times\gB}_{n\ \text{times}}\lra\gB\]

A \emph{formal multilinear function series} over $\gB$ is a sequence
$\omega=(\omega_0,\omega_1,\dots)$, where
$\omega_0\in\widetilde{\gB}$ and $\omega_n\in\mathcal{L}_n(\gB)$ for
$n\geq1$. According to \cite{dykema}, the set of all multilinear
function series over $\gB$ will de denoted by $\mul$.

For $F,G\in\mul$, the \emph{sum} $F+G$ and the \emph{formal product}
$FG$ are the elements from $\mul$ defined by:

\begin{eqnarray*}
(F+G)_n(b_1,\dots,b_n)&=&F_n(b_1,\dots,b_n)+G_n(b_1,\dots,b_n)\\
(FG)_n(b_1,\dots,b_n)&=&
\sum_{k=0}^nF_k(b_1,\dots,b_k)G_{n-k}(b_{k+1},\dots,b_n)\\
\end{eqnarray*}
for any $b_1,\dots,b_n\in\gB$.

If $G_0=0$, then the \emph{formal composition} $F\circ G\in\mul$ is
defined by
\begin{eqnarray*}
(F\circ G)_0&=&F_0\\
 (F\circ
G)_n(b_1,\dots,b_n)&=&\sum_{k=1}^n \sum_{
\substack{{p_1,\dots,p_k\geq1}\\
{p_1+\dots+p_k=n} } } \hspace{-.5cm}
F_k(G_{p_1}(b_1,\dots,b_{p_1}),\dots,\\
&&\hspace{3.7cm}G_{p_k}(b_{q_k+1},\dots,b_{q_k+p_k}))
\end{eqnarray*}
\noindent where $q_j=p_1+\dots+p_{j-1}$, $n\geq1$.

 With the above operations, $\mul$ is an algebra, with
 additional properties similar to ones of power series
 (see \cite{dykema}, \textbf{Proposition 2.3} and
\textbf{Proposition 2.6}):

\begin{prop}\label{mulprop}Let $E,F,G\in\mul$. Then:
\begin{enumerate}
\item[(i)] $1=(1,0,0,\dots)\in\mul$ is a multiplicative
identity element;
\item[(ii)]$F=(F_0,F_1,\dots)$ has a multiplicative inverse
if and only if $F_0$ is an\\
 invertible element of
$\widetilde{\gB}$;
\item[(iii)] if $F_0=0$ and $G_0=0$,
then $(E\circ F)\circ G=E\circ (F\circ G)$;
\item[(iv)]if $G_0=0$, then $(E+F)\circ G=E\circ G+F\circ G$ and
$(EF)\circ G=(E\circ G)(F\circ G)$;
\item[(v)]$I=(0, id_\gB,0,0,\dots)\in\mul$ is an identity element
for the formal composition
\item[(vi)]$F=(0,F_1,F_2,\dots)\in\mul$ has a compositional inverse,
 denoted $F^{\inv}$, if and only if $F_1$ is an invertible element
of $\mathcal{L}_1(\gB)$.
\item[(vii)]if $F=(0,F_1,F_2,\dots)\in\mul$, then
\[(1-F)^{-1}=1+\sum_{k=1}^\infty F^k\]
\end{enumerate}
\end{prop}

For the next definitions and results, $\gA$ will be a $\gB$-algebra
($\gB$  and $\gA$ are not necessarily $\ast$-algebras).

\begin{defn}\emph{For $X\in\gA$ consider}
$\gH_X=(\gH_{X,0}, \gH_{X,1},\dots)\in\mul$, \emph{where}
\begin{eqnarray*}
\gH_{X,0}&=&0\\
\gH_{X,1}(b)&=&b\\
\gH_{X,n}(b_1,\dots,b_n)&=&\Phi(b_1Xb_2\cdots b_{n-1}Xb_n)
\end{eqnarray*}
\emph{for all} $b, b_1,\dots, b_n\in\gB$, $n\geq1$.
\end{defn}

\begin{thm}\label{multhm1}If $X,Y\in\gA$ are monotonically independent over $\gB$,
then
\[\gH_{X+Y}=\gH_X\circ\gH_Y\]
\end{thm}

\begin{proof}
It suffices to show that $\gH_{X+Y,n}=(\gH_X\circ\gH_Y)_n$ for all
$n\geq0$.

For $n=0,1$, the assertion is trivial. For $n\geq2$,
\begin{eqnarray*}
(\gH_X\circ \gH_Y)_n(b_1,\dots,b_n) &=& \sum_{k=1}^n
 \sum
\gH_{X,k}(\gH_{Y,p_1}(b_1,\dots,b_{p_1}),
\dots,\gH_{Y,p_k}(b_{q_k+1},\dots,
b_n))\\
&&\text{where}\ q_j=p_1+\dots+p_{j-1}\
\text{and the second summation}\\
&&\text{ is over all}\ p_1,\dots,p_k\geq1\
\text{such that}\ p_1+\dots+p_k=n\  \\
&=&\Phi\left(\sum_{k=1}^n\sum\gH_{Y,p_1}(b_1,\dots,b_{p_1})X \dots
X\gH_{Y,p_k}(b_{q_k+1},\dots,
b_n))\right)\\
&=&\Phi\left(\sum_{k=1}^n\sum\Phi(b_1Y\dots Yb_{p_1})X\dots
X\Phi(b_{q_k+1}Y\dots Y b_n)\right)
\end{eqnarray*}
Using Definition \ref{defmonindep}, the above relation becomes
\begin{eqnarray*}
(\gH_X\circ \gH_Y)_n(b_1,\dots,b_n) &=&\sum_{k=1}^n\sum
\Phi\left(b_1Y\dots Yb_{p_1}X\dots Xb_{q_k+1}Y\dots Y b_n\right)\\
&&\text{ with the convention that if}\ p_j=1,\ \\
&&\text{then}\ b_{q_j+1}Y\dots Yb_{q_{j+p_j}}=b_{q_j+1}\\
&=&\sum_{\substack{{X_i\in\{X,Y\}}\\{i=1,\dots,n-1}}}
\Phi(b_1X_1\dots X_{n-1}b_n)\\
&=&\Phi\left(b_1(X+Y)\dots
(X+Y)b_n\right)\\
\end{eqnarray*}

On the other hand,
\begin{eqnarray*}
\gH_{X+Y,n}(b_1,\dots,b_n)&=&\Phi\left(b_1(X+Y)\dots
(X+Y)b_n\right)\\
\end{eqnarray*}
hence the conclusion.
\end{proof}

If the algebra $\gA$ is unital, there also exist multilinear
function series analogous  to $\kappa$, $\rho$. First, for
$X\in\gA$, define the elements $\beta_X$ and $\gamma_X$ of $\mul$ by

\begin{eqnarray*}
\beta_{X,0}&=&0\\
\beta_{X,n}(b_1,\dots,b_n)&=&\Phi(b_1Xb_2\dots b_nX)\\
\gamma_{X,0}&=&0\\
\gamma_{X,n}(b_1,\dots,b_n)&=&\Phi(Xb_1X\dots X b_n)\\
\end{eqnarray*}

From the property \ref{mulprop}(ii), the multilinear
 function series
$\gK_X$ and $\mathfrak{r}_X$ are well-defined, where

\begin{eqnarray*}
\gK_X=(1+\beta_X)^{-1}\beta_X\\
\mathfrak{r}_X=\gamma_X(1+\gamma_X)^{-1}.
\end{eqnarray*}

\begin{thm}
Let $U,V\in\gA$ be such that $U-1$ and $V$ are monotonically
independent over $\gB$. Then:
\begin{eqnarray*}
\gK_{VU}=\gK_U\circ\gK_V\\
\mathfrak{r}_{UV}=\mathfrak{r}_V\circ\mathfrak{r}_U.
\end{eqnarray*}
\end{thm}

The proof is a routine (though tedious) verification, using
Proposition \ref{mulprop} and the techniques from the proof of
Theorem \ref{multhm1} and  Theorem \ref{thm2}.

\section{Semi-Inner Product Bimodules}\label{fock}

 The terminology in used in this section is the one from \cite{Lance}.
 Let $\gB$ be a unital $C^*$-algebra. A \emph{semi-inner product $\gB$-bimodule}
 is a linear space $\cE$ which is a $\gB$-bimodule, together with a
 map $(x,y)\mapsto \langle x,y \rangle:\cE\times\cE\lra \gB$ such
 that:

\begin{enumerate}
\item[(i)] $\langle x, \alpha y +\beta z\rangle=
\alpha\langle x,y\rangle+ \beta\langle x,z\rangle$ for any
$x,y,z\in\cE, \alpha,\beta\in\mathbb{C}$.
\item[(ii)]$\langle x,ya\rangle =\langle x,y\rangle a$,\
for any $x,y\in\cE, a\in\gB$.
\item[(iii)]$\langle y , x\rangle=\langle x, y\rangle^*$\
for any $x,y\in\cE$
\item[(iv)]$\langle x,x\rangle\geq 0$\ for any $x\in\cE$.
\end{enumerate}

 $\cE$ is called an
\emph{inner-product $\gB$-bimodule} if $\langle x,x\rangle =0$
implies $x=0$ and \emph{Hilbert $\gB$-bimodule} if it is complete
with respect to the norm $||x||_0=||\langle
x,x\rangle||^\frac{1}{2}$. (second norm is the C$^*$-algebra norm of
$\gB$.)
 The algebra of $\gB$-linear (not necessarily bounded) operators on
 $\cE$ will be denoted by $\cL(\cE)$.

 \subsection{}Given a family
 $(\cE_i)_{i\in \I}$  of semi-inner product $\gB$-bimodules
 indexed by a totally ordered set
 $\I\subseteq \mathbb{Z}$,
 we define, following \cite{muraki1} and \cite{pimsner},
 the monotonic product $\mE$ of $(\cE_i)_{i\in \I}$
 to be the semi-inner product $\gB$-bimodule:
 \[\mE=\gB\oplus
\left(\bigoplus_{n\geq1}
 \bigoplus_{(i_1,\dots, i_n)\in \orIn}
\cE_{i_1}\otimes\cE_{i_2}\otimes\dots\otimes\cE_{i_n}
 \right)\]

 where all the tensor products are with amalgamation over $\gB$ and
 \[\orIn=\{(i_1,\dots,i_n):\ i_1, \dots i_n\in \I, i_1>\dots>i_n\},\]

with the inner-product given by
\[\langle f_1\otimes\dots\otimes f_n, e_1\otimes\dots \otimes e_m\rangle
=\delta_{m,n}\langle f_n, \langle f_{n-1}\dots,\langle f_1,
e_1\rangle\dots e_{n-1}\rangle e_n\rangle.\]

Note that in general $\cE^m$ is not an inner-product $\gB$-bimodule
even if $\cE_i$ are inner-product bimodules or Hilbert bimodules.
For example, if $\langle f_1,f_1 \rangle = b^*b>0$ and $bf_2=0$,
then
\[\langle f_1\otimes f_2, f_1\otimes f_2 \rangle =
\langle f_2, \langle f_1,f_1\rangle f_2\rangle=\langle f_2, b^*b
f_2\rangle =0;\]
see also \cite{speicher1}.

 If $i\in \I$ is fixed, we have the natural identification
 \[\mE=\left(\left(\gB\oplus\cE_i\right)\otimes_{\gB}
 \left(\gB\oplus\mE_{(<i)}\right)\right)\oplus\mE_{(>i)}\]
 where
 \begin{eqnarray*}
 \mE_{(<i)}&=&\bigoplus_{n\geq1}
 \bigoplus_{
 \substack
 {
 (i_1,\dots, i_n)\in \orIn\\
 i_1<i
 }
 }
\cE_{i_1}\otimes\cE_{i_2}\otimes\dots\otimes\cE_{i_n}\\
 \mE_{(>i)}&=&\bigoplus_{n\geq1}
 \bigoplus_{
 \substack
 {
 (i_1,\dots, i_n)\in \orIn\\
 i_1>i
 }
 }
\cE_{i_1}\otimes\cE_{i_2}\otimes\dots\otimes\cE_{i_n}.
 \end{eqnarray*}

Based on this decomposition, one also has the (non-unital)
$\ast$-representation\\ $\lambda_i:\cL(\gB\oplus\cE_i)\lra\cL(\mE)$
\[\lambda_i(A)=\left(A\otimes I_{\gB\oplus\mE_{(<i)}}\right)
\oplus 0_{\mE_{(>i)}}\]

\begin{thm}With the above notations,
$\{\lambda_i(\cL(\gB\oplus\cE_i))\}_{i\in\I}$
 are monotonically independent in $\cL(\mE)$
with respect to the conditional expectation $\Phi(\cdot)=\langle
1,\cdot 1\rangle$.
\end{thm}
\begin{proof}
We need to show that the two conditions from the definition of
monotonic independence (Definition \ref{defmonindep}) are satisfied.
In fact, it will be shown that the family
$\{\lambda_i(\cL(\gB\oplus\cE_i))\}_{i\in\I}$ satisfies
\ref{defmonindep}(b) and the stricter condition
\ref{defmonindep}(a$^\prime$).

The proof is similar to the proof of Theorem  2.1  from
\cite{muraki1}. For $i\in\I$, consider $A_i\in\cL(\gB\oplus\cE_i)$
and $X_i=\lambda_i(A_i)$.

We can write:
\[X_i1=\alpha_i+ s_i,\ \alpha_i\in\gB, s_i\in\cE_i\\\]
\noindent where $1\in\gB\subset\gB\oplus\cE_{i,j}$.

 If $k<l$,
\begin{eqnarray*}
X_kX_l1&=& X_k(\alpha_l+ s_l)\\
&=&X_k1\alpha_l\\
&=&X_k\langle 1,X_l 1\rangle
\end{eqnarray*}
therefore

\[X_jX_{k_1}\dots X_{k_n}1
=\langle 1,X_j 1\rangle \langle 1,X_{k_1}1\rangle\dots
\langle1,X_{k_n}1\rangle\]

\noindent whenever $j<k_1<\dots <k_n$.

Also, writing $\mE=\gB\oplus \cE^0$, note that $X_lf\in\cE^0$ for
any $f\in \cE^0$ and any $l\in\I$, and that ($k<l$):
\begin{eqnarray*}
X_lX_k1&=& X_l(\alpha_k + s_k)\\
&=&X_l1\alpha_k+X_l(1\otimes s_k)\\
&=&X_l\langle 1,X_k 1\rangle+ (\alpha_l+s_l)\otimes s_k\\
&=& X_l\langle 1,X_k 1\rangle+ f\ \ \text{for some} f\in\cE^0.
\end{eqnarray*}

Iterating the above relations, for $i_m>\dots>i_1>j<k_1<\dots<k_n$,
we obtain

\begin{eqnarray*}
\langle1, X_{i_m}\dots X_{i_1}X_j X_{k_1}\dots X_{k_n}1\rangle
&=&\langle1, X_{i_m}\dots X_{i_1}1\rangle\langle1,X_j1\rangle
\dots\langle1,X_{k_n}1\rangle\\
&=&\left(\langle1,X_{i_m}1\rangle\dots\langle1,X_{i_1}1
\rangle+\langle1,f\rangle\right)
\langle1,X_j1\rangle
\dots\langle1,X_{k_n}1\rangle\\
&=&\langle1,X_{i_m}1\rangle\dots\langle1,X_{k_n}1\rangle
\end{eqnarray*}
\noindent that is, property (b).

For  $i<j>k$, a direct computation gives

\begin{eqnarray*}
 X_iX_jX_k1&=&X_iX_j(\alpha_k+s_k)\\
&=&X_i(X_j1)\alpha_l+X_iX_j(1\otimes s_k)\\
&=&X_i(\alpha_j+s_j)\alpha_k+X_i(\alpha_j+s_j)\otimes s_k\\
&=&X_i(\alpha_j\alpha_k)+X_i(\alpha_js_k)+X_i(s_j\otimes s_k)\\
&=&X_i\alpha_j(\alpha_k+s_k)\\
&=&X_i\langle1,X_j1\rangle X_k1,
\end{eqnarray*}

\noindent so it remains to show (a$^\prime$) on elements of the form
$\widetilde{h}=h_{i_1}\otimes\dots\otimes h_{i_n}$,
$h_{i_l}\in\cE_{i_l}$.

If $i_1>i$, then $X_2h_{i_1}\otimes\dots\otimes h_{i_n}=0$,
therefore
\[
X_1YX_2\widetilde{h}=0
=X_1\langle1,Y1\rangle X_2\widetilde{h}.
\]

If $i_1=i$, with the notations $h^0=h_{i_2}\otimes\dots\otimes
h_{i_n}$ and $X_2h_{i_1}=\theta\oplus u$ for some
 $\theta\in\gB,u\in\cE_i$,
 one has
 \begin{eqnarray*}
 X_1YX_2h_{i_1}\otimes\dots\otimes h_{i_n}
 &=&
 X_1Y(\theta\oplus u)h^0\\
 &=&
 X_1[\beta\theta\oplus t \theta \oplus
 (\beta\oplus t)\otimes u]\otimes h^0\\
 &=&
 X_1[\beta\theta+\beta u]\otimes h^0\\
 &=&
 X_1 \beta (\theta\oplus u)\otimes h^0\\
 &=&
 X_1\langle1,Y1\rangle X_2h_{i_1}\otimes\dots\otimes h_{i_n}.\\
 \end{eqnarray*}
 The case $i_1<i$ is analogous.

\end{proof}

\subsection{}\label{weaklymonotone}

The \emph{weakly monotone product} of the bimodules
$\{\cE_i\}_{i\in\I}$ is the semi-inner product $\gB$-bimodule

\[\wme=\gB\oplus\bigoplus_{n=1}^\infty
\left(\bigoplus_{i_1\geq\dots\geq i_n}
\cE_{i_1}\otimes\dots\otimes\cE_{i_n}\right)\]

If $\I$ has only one element, $i_0$, then $\wme$ is the full Fock
bimodule over $\cE_{i_0}$, $\mathcal{F}(\cE_{i_0})$ (see
\cite{pimsner}, \cite{speicher1}).

For $j\in\I$, let $\mathfrak{J}=\{l\in\I, l\leq j\}$ and let
$\wme(j)$ be the weakly monotonic product of
$\{\cE_l\}_{l\in\mathfrak{J}}$. We will also use the following
notations
\begin{eqnarray*}
\mathcal{F}_0(\cE)&=&\mathcal{F}(\cE)\ominus \gB\\
\wme_0(\cE)&=&\wme\ominus \gB\\
\wme_0(j)&=&\wme(j)\ominus \gB
\end{eqnarray*}

 For
$f\in\cE_i$, define the $\gB$-linear creation and annihilation maps
$a^*(f)$ and $a(f)$ on $\wme$ by:

\begin{eqnarray*}
a^*(f)1&=&f\\
a^*(f)f_{i_1}\otimes\dots\otimes f_{i_n}
 &=&
 \left\{
 \begin{matrix}
 f\otimes f_{i_1}\otimes\dots\otimes f_{i_n}, & \text{if}\ i\geq i_1\\
 0, & \text{if}\ i<i_1\\
 \end{matrix}
 \right.\\
 a(f)1&=&0\\
 a(f)f_{i_1}\otimes\dots\otimes f_{i_n}
 &=&
 \left\{
 \begin{matrix}
 \langle f, f_{i_1}\rangle f_{i_2}\otimes\dots\otimes f_{i_n}, & \text{if}\ i=i_1\\
 0, & \text{if}\ i\neq i_1\\
 \end{matrix}
 \right.\\
\end{eqnarray*}

Note that $a(f)$ and $a^*(f)$ are adjoint to each other. Denote by
$G(f)$ their sum, $G(f)=a(f)+a^*(f)$, and by $\gA_i$ the algebra
generated over $\gB$ by $\{G(f),\ f\in\cE_i\}$.

 We will use the shorthand notation
  $\Phi(\cdot)$ for the $\gB$-valued functional $\langle1,\cdot1\rangle$ on the set of all
  $\gB$-linear maps  on $\wme$.
  Also, for
  $\widetilde{e}=e_1\otimes\dots\otimes e_n, \
  (e_l\in\cE_k, 1\leq l \leq n)$
  we will use the notations
  \begin{eqnarray*}
  A^*(\widetilde{e})&=&a^*(e_1)\cdots a^*(e_n)\\
  A(\widetilde{e})&=&a(e_1)\cdots a(e_n)
  \end{eqnarray*}

\begin{lemma}\label{lemmacreannih}
For any $f_1,\dots,f_n\in\cE_k$ there are some sequences  of
elements of $\wme_0(k)$,
$(\widetilde{e}_{r})_{r=1}^{N_1},(\widetilde{g}_s)_{s=1}^{N_2},
(\widetilde{h}_q)_{q=1}^{N_3}, (\widetilde{k}_q)_{q=1}^{N_3}$ , such
that
\[P=\prod_{l=1}^nG(f_l)\]
can be written as:
\begin{equation}\label{creannih}
P=\Phi(P) +\sum_{r=1}^{N_1}A^*(\widetilde{e}_r)
+\sum_{s=1}^{N_2}A(\widetilde{g}_s)
+\sum_{q=1}^{N_3}A^*(\widetilde{h}_q)A(\widetilde{k}_q)
\end{equation}
\end{lemma}

\begin{proof}
Let (\ref{creannih})$^\prime$ be the weaker form of (\ref{creannih})
where $\Phi(P)$ is replaced by some element  $\alpha\in\gB$. Note
that (\ref{creannih})$^\prime$ is in fact equivalent to
(\ref{creannih}), since

\begin{eqnarray*}
\Phi(P)&=& \langle1,P1\rangle\\
&=&\langle1,\alpha1 +\sum_{r=1}^{N_1}A^*(\widetilde{e}_r)1
+\sum_{s=1}^{N_2}A(\widetilde{g}_s)1
+\sum_{q=1}^{N_3}A^*(\widetilde{h}_q)A(\widetilde{k}_q)1\rangle\\
&=& \langle1,1\alpha+ \sum_{r=1}^{N_1}\widetilde{e}_r\rangle\\
&=&\alpha
\end{eqnarray*}

It remains to prove (\ref{creannih})$^\prime$.

Note that
\begin{eqnarray*}
P&=&\prod_{l=1}^nG(f_l)
=\prod_{l=1}^n(a^*(f_l)+a(f_l))\\
&=&\sum_{\begin{array} {clcr}(\varepsilon_1,\dots,\varepsilon_n)\\
\varepsilon_l\in\{1,2\}
 \end{array}}
a_{\varepsilon_1}(f_1)\dots a_{\varepsilon_n}(f_n)
\end{eqnarray*}

where $a_1$ stands for $a$ and $a_2$ stands for $a^*$.

Also, for any $f,g,h\in\cE_k, \alpha\in\gB$ and
$\varepsilon\in\{1,2\}$
\begin{eqnarray*}
a(f)a^*(g)&=&\langle f,g\rangle I\\
a_\varepsilon(h)\langle f,g\rangle&=&a_\varepsilon(h\langle
f,g\rangle)\\
\alpha a(f)&=&a(\alpha^*f)\\
\alpha a^*(f)&=&a^*(\alpha f)
\end{eqnarray*}

It follows that in the expression of $a_{\varepsilon_1}(f_1)\dots
a_{\varepsilon_n}(f_n)$ any $a(f_p)a^*(f_{p+1})$ can be reduced to
$\langle f_p,f_{p+1}\rangle$ which can be included in the expression
of the previous or following factor. Iterating, after a finite
number of steps no summand will have factors of the type $a(f_q)$ in
front of factors of the type $a^*(f_p)$, so
(\ref{creannih})$^\prime$ is proved.

\end{proof}

\begin{lemma}
\label{formajut} Any $X\in\gA_i$ is satisfying the following
properties:
\begin{enumerate}
\item[(i)]$\wme(i)$ is $X$-invariant
\item[(ii)]$X1=\Phi(X)+s$, for some $s\in\mathcal{F}_0(\cE_i)$
\item[(iii)]if $u\in\wme(j), j<i,$ then
\[Xu=\Phi(X)u+t\otimes u, \ \text{for some}\ t\in\mathcal{F}_0(\cE_i)\]
\item[(iv)]if $i<j$ and  $v\in\mathcal{F}_0(\cE_j)\otimes\wme(k),
k<j,$ then
\[Xv=0 \hspace{5cm}\]
\end{enumerate}
\end{lemma}

\begin{proof}

\noindent\\

 \textbf{(i)}It is enough to verify that the property holds true for $X=G(f)$, $f\in\cE_i$.
 Indeed,
\[G(f)1=f\in\mathcal{F}(\cE_i)\subset\wme(i)\]
 and for any $f_1\otimes\dots\otimes f_n\in\wme(i)$,
\[
G(f)f_1\otimes\dots\otimes f_n=f\otimes f_1\otimes\dots\otimes f_n
+\langle f,f_1\rangle f_2\otimes\dots\otimes f_n\in\wme(i).
\]

\textbf{(ii)}If $f,f_1,\dots,f_n\in\mathcal{F}(\cE_i)$

\[
G(f)f_1\otimes\dots\otimes f_n=f\otimes f_1\otimes\dots\otimes f_n
+\langle f,f_1\rangle f_2\otimes\dots\otimes
f_n\in\mathcal{F}(\cE_i).
\]
It follows that $\mathcal{F}(\cE_i)$ is invariant to $\gA_i$. Since
$1\in\mathcal{F}(\cE_i)$, we have $X1\in\mathcal{F}(\cE_i)$, and the
conclusion follows from the orthogonality of $\gB$ and
$\mathcal{F}_0(\cE_i)$.

\textbf{(iii)} It is enough to prove the relation for $X=G(f_1)\dots
G(f_n)$, $f_l\in\cE_i$. First note that for any
$\widetilde{f}\in\mathcal{F}_0(\cE_i)$

\begin{eqnarray*}
A^*(\widetilde{f})u&=&\widetilde{f}\otimes u \\
A(\widetilde{f})u&=&0
\end{eqnarray*}

 From Lemma \ref{lemmacreannih}, there are some
sequences
$(\widetilde{e}_{r})_{r=1}^{N_1},(\widetilde{g}_s)_{s=1}^{N_2},
(\widetilde{h}_q)_{q=1}^{N_3}, (\widetilde{k}_q)_{q=1}^{N_3}$, of
elements of $\mathcal{F}(\cE_i)$ , such that

 \begin{eqnarray*}
 Xu
 &=&
 \Phi(X)u +\sum_{r=1}^{N_1}A^*(\widetilde{e}_r)u
+\sum_{s=1}^{N_2}A(\widetilde{g}_s)u
+\sum_{q=1}^{N_3}A^*(\widetilde{h}_q)A(\widetilde{k}_q)u\\
&=& \Phi(X)u +\sum_{r=1}^{N_1}\widetilde{e}_r\otimes u\\
&=& \Phi(X)u +\left(\sum_{r=1}^{N_1}\widetilde{e}_r\right)\otimes u
 \end{eqnarray*}

\textbf{(iv)}Similarly, it is enough to prove the relation for
$X=G(f_1)\dots G(f_n)$, $f_l\in\cE_i$. First note that for any
$\widetilde{f}\in\mathcal{F}_0(\cE_i)$

\[
A^*(\widetilde{f})v= A(\widetilde{f})v=0
\]
 From Lemma \ref{lemmacreannih}, there are some
sequences
$(\widetilde{e}_{r})_{r=1}^{N_1},(\widetilde{g}_s)_{s=1}^{N_2},
(\widetilde{h}_q)_{q=1}^{N_3}, (\widetilde{k}_q)_{q=1}^{N_3}$, of
elements of $\mathcal{F}(\cE_i)$ , such that

 \begin{eqnarray*}
Xv
 &=&
 \Phi(X)v +\sum_{r=1}^{N_1}A^*(\widetilde{e}_r)u
+\sum_{s=1}^{N_2}A(\widetilde{g}_s)v
+\sum_{q=1}^{N_3}A^*(\widetilde{h}_q)A(\widetilde{k}_q)v\\
&=& \Phi(X)v
 \end{eqnarray*}

\end{proof}

\begin{thm}
The algebras $\{\gA_i\}_{i\in\I}$  are monotonically independent
with respect to the $\gB$-valued functional
$\Phi(\cdot)=\langle1,\cdot 1\rangle$.
\end{thm}

\begin{proof}
Let $X_i\in\gA_i, i\in\I$. We will prove that they satisfy the
relations (b) and (a$^\prime$) from  the definition of the monotonic
independence. If $k<l$, from Lemma \ref{formajut}, we have
\begin{eqnarray*}
X_kX_l1&=&X_k(\Phi(X_l)+t) \ \text{for some}\
t\in\mathcal{F}_0(\cE_l)\\
&=&X_k1\Phi(X_l), \ \text {since}\ Xt=0,
\end{eqnarray*}
therefore
\[X_jX_{i_1}\dots X_{k_n}1=X_j1\Phi(X_{k_1})\dots\Phi(X_{k_n})\]
\noindent whenever $j<k_1<\dots k_n$. Similarly, ($k<l$):

\begin{eqnarray*}
X_lX_k1&=&X_l(\Phi(X_k)+t_k,\ \text{for some}\
t_k\in\mathcal{F}_0(\cE_k)\\
&=&\Phi(X_l)\Phi(X_k)+\Phi(X_l)t_k+t_l\otimes(\Phi(X_k)+t_k)\
\text{for some}\ t_l\in\mathcal{F}_0(\cE_l)\\
&=&\Phi(X_l)\Phi(X_k)+t,\ \text{for some}\ t\in\wme_0(l).
\end{eqnarray*}

Using the above relations, we obtain
\begin{eqnarray*}
\Phi(X_{i_1}\dots X_{1_n}X_jX_{k_1}\dots X_{k_m})
&=&
\Phi(X_{i_1}\dots X_{1_n}X_j)\Phi(X_{k_1})\dots \Phi(X_{k_m})\\
&=&\Phi(\Phi(X_{i_1})\dots \Phi(X_j)+s)\Phi(X_{k_1})\dots
\Phi(X_{k_m},)\\
&&\hspace{3cm} \text{for some}\ s\in\wme_0(k_m)\\
&=&\Phi(X_{i_1})\dots \Phi(X_j)+s)\Phi(X_{k_1})\dots \Phi(X_{k_m}),
\end{eqnarray*}
\noindent that is, property (b).

Also, for $i<j>k$,

\begin{eqnarray*}
X_iX_jX_k1&=&X_iX_j\left(\Phi(X_k)+s\right)\ \text{for some}\
s\in\mathcal{F}_0(\cE_k)\\
&=&(X_iX_j1)\Phi(X_k)+X_iX_js)\\
&=&X_i\left(\Phi(X_j)+t_1\right)\Phi(X_k)+X_i\left(\Phi(X_j)s+t_2\otimes
s\right),\ \text{for some}\ t_1,t_2\in\mathcal{F}_0(\cE_j)\\
&=&X_i\Phi(X_j)\Phi(X_k)+X_i\Phi(X_j)s\\
&=&X_i\Phi(X_j)\left(\Phi(X_k)+s\right)\\
&=&X_i\Phi(X_j)X_k1.
\end{eqnarray*}

It remains to prove that $X_iX_jX_k$ and $X_i\Phi(X_j)X_k$ coincide
on vectors of the form $\widetilde{f}=f_{i_1}\otimes\dots\otimes
f_{i_m}$, where $f_{i_k}\in\cE_{i_k}, i_1\geq\dots\geq i_m$.

If $i_1>k$, then $X_k\widetilde{f}=0$, so the equality is trivial.

If $i_1\leq k$, then Lemma \ref{formajut} implies that
$X_k\widetilde{f}\in\wme(k)$, therefore:
\begin{eqnarray*}
X_iX_jX_k\widetilde{f}&=&X_i\left(\Phi(X_j)X_k\widetilde{f}+t\otimes
X_k\widetilde{f} \right),\ \text{for some}\ t\in\mathcal{F}_0(\cE_j)\\
&=&X_i\Phi(X_j)X_k\widetilde{f}
\end{eqnarray*}
\end{proof}

\begin{remark}\emph{An analogous construction can be done for creation and
annihilation maps on $\cE^m$, and similar computations will lead to
the monotonic independence of the correspondent algebras (see
\cite{muraki2}).}
\end{remark}

\section{central limit theorem}

  In this section $\gA$ will be a $\ast$-algebra,
   $\gB$ a subalgebra of $\gA$ which is also a $C^*$-algebra
   and $\Phi:\gA\lra\gB$ a conditional
  expectation. $\gB_+\langle\xi\rangle$ will be the $\ast$-algebra
  generated by $\gB$ and the selfadjoint symbol $\xi$, as described
  in Introduction.

As discussed in Section 2, given $X$ a selfadjoint element of $\gA$,
the $n$-th moment of $X$ is the multilinear function
$m_{X,n}:\gB^{n-1}\lra\gB$,
\[
 m_{X,n}(b_1,\dots,b_{n-1})=\Phi(Xb_1X\cdots Xb_{n-1}X).
\]

\noindent We define \emph{the moment function} of $X$ as
\[\mu_X=\bigoplus_{m=1}^\infty\mu_X^{(m)}\]

Before stating the main theorem of this section, we will begin with
some combinatorial considerations on the joint moments of the family
of selfadjoint elements $(X_n)_{n\geq 1}$ of $\gA$ with the
properties:

\begin{enumerate}
\item[(1)] for any $i<j$, $X_i$ and $X_j$ are monotonically
independent over $\gB$;
\item[(2)] all $X_k$ have the same moment function, denoted by $\mu$.
\end{enumerate}

 Let $NC(m)$ be the set of all non-crossing partitions of the
 ordered set $\{1,2,\dots,m\}$. For $\g\in NC(m)$, let
 $B=\{b_1,b_2,\dots,b_p\}$ and $C=\{c_1,c_2,\dots,c_q\}$ be two
 blocks of $\g$. We say that $C$ is \emph{interior to} $B$ if there
 is a index $k\in\{1,\dots, p-1\}$ such that
 $b_k<c_1,c_2,\dots,c_q<b_{k+1}$. $B$ and $C$ will be called \emph{adjacent}
  if $c_1=b_p+1$ or $b_1=c_q+1$. The block $B$ will be called \emph{outer}
  if it is not interior to any other block of $\g$.

 To each $m$-tuple $(i_1,\dots,i_m)$ of indices from $\{1,2,,\dots\}$ we associate a
 unique non-crossing partition $nc[i_1,\dots,i_m]\in NC(m)$ as follows:
 \begin{enumerate}
 \item[(1)]if $m=1$, then $nc[i_1]=(1)$
 \item[(2)] if $m>1$, put
 \[B=\{k,
 i_k=min\{i_1,\dots,i_m\}\}=\{k_1,\dots,k_p\}\]
\noindent  and define
 \[nc[i_1,\dots,i_m]=B\sqcup nc[i_1,\dots,i_{k_1}-1]\sqcup
 nc[i_{k_1}+1,\dots,i_{k_2}-1]\sqcup\dots
 \sqcup nc[i_{k_p}+1,\dots,i_m]\]
 \end{enumerate}

Reciprocally, the $m$-tuple $(i_1,\dots,i_m)$ will be called an
admissible configuration for $\g\in NC(m)$ if
$nc[i_1,\dots,i_m]=\g$.

\begin{lemma}\label{5lemma1}
Suppose $(i_1,\dots, i_m)$ is an admissible configuration for $\g\in
NC(m)$ and $B=\{k_1,\dots,k_p\}$ is an outer block of $\g$. Then,
for any $b_1,\dots,b_{m-1}\in\gB$, $mu$ the common moment function
of $\{X_n\}_n$, we have

\begin{eqnarray*}
\Phi(X_{i_1}b_1 X_{i_2}\dots b_{m-1} X_{i_m})&&\\
&&\hspace{-3cm}=
 \mu(\Phi(X_{i_1}b_1\dots
X_{i_{k_1-1}}b_{k_1-1})
 ,\Phi(b_{k_1}X_{i_{k_1+1}}\dots b_{k_2}),\dots,
 \Phi(b_{k_p}\dots X_{i_m}))\\
\end{eqnarray*}

\end{lemma}

\begin{proof}
If $\g$ has only one block, then the result is trivial. If $\g$ has
more than one block, but only one outer block, $B$, then $B=\{k:
i_k=min\{i_1,\dots,i_m\}\}=\{k_1,\dots, k_p\}$, since the last set
forms always an outer block. Also this block must contain $1$ and
$m$. The monotonic independence of $(X_n)_{n\geq1}$ over $\gB$
implies
\begin{eqnarray*}
\Phi(X_{i_1}b_1 X_{i_2}\dots b_{m-1} X_{i_m})&&\\
&&\hspace{-3.4cm}=
 \Phi(X_{i_1}\Phi(b_1X_{i_2}b_2\dots
X_{i_{{k_2-1}}}b_{k_2-1})X_{i_{k_2}}
 \Phi(b_{k_2}X_{i_{k_2+1}}\dots)\dots
 X_{i_m})\\
&&\hspace{-3.4cm}= \mu(\Phi(b_1X_{i_2}\dots b_{k_2-1}),\dots,
\Phi(b_{k_{p-1}}X_{i_{k_{p-1}+1}}\dots b_{m-1}))
\end{eqnarray*}

If $\g$ has more than one outer block, the result comes by induction
on the number of blocks of $\g$. Suppose the result is true for less
than $r$ blocks and that $\g$ has exactly $r$ blocks. Consider again
$B_0=\{k, i_k=min\{i_1,\dots,i_m\}\}=\{k_1,\dots, k_p\}$. Using
again Definition \ref{defmonindep}, we obtain

\begin{eqnarray*}
\Phi(X_{i_1}b_1 X_{i_2}\dots b_{m-1} X_{i_m})&&\\
&&\hspace{-3.4cm}= \Phi(\Phi(X_{i_1}b_1\dots
b_{k_1-1})X_{i_{k_1}}\dots X_{i_{k_m}}
\Phi(b_{k_m}X_{i_{k_m+1}}\dots X_{i_m}))\\
&&\hspace{-3.4cm}=
 \mu(\Phi(X_{i_1}b_1\dots
b_{k_1-1}),\dots ,\Phi(b_{k_m}X_{i_{k_m+1}}\dots X_{i_m}))
\end{eqnarray*}

If $B=B_0$, then the result is proved above. If $B\neq B_0$, then
without losing generality we can suppose that $B$ is at the right of
$B_0$, hence

\begin{eqnarray*}
\Phi(X_{i_1}b_1 X_{i_2}\dots b_{m-1} X_{i_m})&&\\
&&\hspace{-3.4cm}= \Phi(\Phi(X_{i_1}b_1\dots
b_{k_1-1})X_{i_{k_1}}\dots X_{i_{k_m}}
\Phi(b_{k_m}X_{i_{k_m+1}}\dots X_{i_m}))\\
&&\hspace{-3.4cm}= \Phi(X_{i_1}b_1\dots X_{i_{k_1-1}})
\Phi(b_{k_1-1}X_{i_{k_1}}\dots X_{i_m}),
\end{eqnarray*}

\noindent and the result follows applying the induction hypothesis
to $\Phi(X_{i_1}b_1\dots X_{i_{k_1-1}})$.

\end{proof}

\begin{lemma}\label{5lemma2}
If $(i_1,\dots,i_m)$ and $(l_1,\dots,l_m)$ are two admissible
configurations for $\g\in NC(m)$, then for any
$b_1,\dots,b_{m-1}\in\gB$, one has:
\[\Phi(X_{i_1}b_1 X_{i_2}\dots b_{m-1} X_{i_m})
=\Phi(X_{l_1}b_1 X_{l_2}\dots b_{m-1} X_{l_m})\]
\end{lemma}

\begin{proof}
 Again, if $\g$ is the partition with a single block, then $i_1=\dots=i_m$
 and $l_1=\dots = l_m$ and

\begin{eqnarray*}
\Phi(X_{i_1}b_1 X_{i_2}\dots b_{m-1} X_{i_m})
 &=&
 \Phi(X_{i_1}b_1 X_{i_1}\dots b_{m-1} X_{i_1})\\
 &=&
 \mu(b_1,\dots, b_{m-1})\\
 &=&
\Phi(X_{l_1}b_1 X_{l_2}\dots b_{m-1} X_{l_m})
\end{eqnarray*}

The conclusion follows now by induction on $m$.

If $m=1$, then $\Phi(X_{i_1})=\Phi(X_{l_1})$.

Suppose the result is true for $m\leq N-1$ and that $\g\in NC(N)$
has more than one block. Let then $B=\{k_1,\dots, k_p\}$ be a outer
block of $\g$. From Lemma \ref{5lemma1}
\begin{eqnarray*}
\Phi(X_{i_1}b_1 X_{i_2}\dots b_{m-1} X_{i_m})&&\\
 &&\hspace{-2.5cm}=
\mu(\Phi(X_{i_1}b_1\dots X_{i_{k_1-1}}b_{k_1-1})
 ,\Phi(b_{k_1}X_{i_{k_1+1}}\dots b_{k_2}),\dots,
 \Phi(b_{k_p}\dots X_{i_m}))\\
&&\hspace{-2.5cm}= \mu(\Phi(X_{l_1} b_1\dots X_{l_{k_1-1}}
b_{k_1-1})
 ,\Phi(b_{k_1}X_{l_{k_1+1}}\dots b_{k_2}),\dots,
 \Phi(b_{k_p}\dots X_{l_m}))\\
 &&\hspace{-2.5cm}=\Phi(X_{l_1}b_1 X_{l_2}\dots b_{m-1} X_{l_m})
\end{eqnarray*}

\end{proof}

Since the value $\Phi((X_{i_1}b_1 X_{i_2}\dots b_{m-1} X_{i_m}))$ is
the same for all the admissible configurations $(i_1,\dots,\i_m)$,
we will denote it by $V(\g,b_1,\dots, b_{m-1})$.

\begin{thm}\label{CLT}
Let  $(X_n)_{n=1}^\infty$ be a sequence of selfadjoint elements from
$\gA$ such that:
\begin{enumerate}
\item[(1)] $\{X_n\}_n$ is a monotonically independent family
\item[(2)]
$\mu_{X_i}=\mu_{X_j}\ \text{for any $i,j\geq1, n\geq0$}$

\item[(3)] $\Phi(X_k)=0$
\end{enumerate}
Then there exists a conditional expectation
\[\nu:\gB_+\langle\xi\rangle\lra\gB\]
with the property:
\begin{equation}
\lim_{N\to \infty}\Phi\left(f\left(\frac{X_1+\dots
X_N}{\sqrt{N}}\right)\right)=\nu(f)
\end{equation}
for any $f\in\gB_+\langle\xi\rangle$.
 Moreover $\nu(f)$ depends only on the second order moments of
 $X_i$.
\end{thm}

\begin{proof}
For convenience, we will use the notations $\mu$ for $\mu_{X_i}$
($i\geq1$), $\adm(\g)$ for the set of all admissible configurations
of $\g$, $\adm(\g,N)$ for the set of all admissible configurations
of $\g$ with indices from $\{1,2,\dots, N\}$, and $PP(m)$ for the
set of all non-crossing pair partitions (partitions where each block
has exactly two elements) of $\{1,\dots,m\}$.

It is enough to show the property for some arbitrary
$b_1,\dots,b_{m-1}\in\gB$ and
\[f=\xi b_1\xi\dots b_{m-1}\xi\]
From Lemma \ref{5lemma2}, one has:

\begin{eqnarray*}
\Phi\left(\left(\frac{X_1+\dots X_N}{\sqrt{N}}\right)b_1\dots
b_{m-1}\left(\frac{X_1+\dots X_N}{\sqrt{N}}\right)\right)\\
&&\hspace{-5cm}=\frac{1}{N^{\m}} \sum_{(i_1,\dots,i_m)}
\Phi(X_{i_1}b_1\dots b_{m-1}X_{i_m})\\
&&\hspace{-5cm}=\frac{1}{N^{\m}}
\sum_{\g\in NC(m)}V(\g,b_1,\dots,b_{m-1}) card(\adm(\g,N))\\
\end{eqnarray*}

 If $\g$ contains blocks with only one element, the condition
 $\Phi(X_i)=0$ ($i\geq1$) and Lemma \ref{5lemma1} imply
 that $V(\g,b_1, \dots, b_{m-1})=0$.

 Also, if $\g$ has less than $\m$ blocks, since
 \[card(\adm(\g,N))<N^{card(\g)}<N^\m\]
 we have that
 \[\lim_{N\to \infty}
 \frac{1}{N^{\m}}V(\g,b_1,\dots,b_{m-1})card(\adm(\g,N))=0.
 \]

 If follows that only the pair partitions contribute to the limit,
 that is

 \begin{eqnarray*}
 \lim_{N\to\infty}
\Phi\left(\left(\frac{X_1+\dots X_N}{\sqrt{N}}\right)b_1\dots
b_{m-1}\left(\frac{X_1+\dots
X_N}{\sqrt{N}}\right)\right)&&\\
&&\hspace{-5.3cm}= \lim_{N\to\infty}\frac{1}{N^{\m}} \sum_{\g\in
PP(m)}V(\g,b_1,\dots,b_{m-1})card(\adm(\g,N))
 \end{eqnarray*}

In particular, for $m$ odd, the limit exists and it is equal to
zero.

If $m$ is even, note first that
\[\adm(\g,N)=\bigsqcup_{k=1}^{\m}\adm(\g,N,k)\]
and that
\[card(\adm(\g,N,k))={{N}\choose{k}}card(\adm(\g,k,k))\]

therefore
\begin{eqnarray*}
\lim_{N\to\infty} \Phi\left(\left(\frac{X_1+\dots
X_N}{\sqrt{N}}\right)b_1\dots b_{m-1}\left(\frac{X_1+\dots
X_N}{\sqrt{N}}\right)\right)&&\\
&&\hspace{-7.3cm}= \lim_{N\to\infty}\frac{1}{N^{\m}} \sum_{\g\in
PP(m)}V(\g,b_1,\dots,b_{m-1})\sum_{k=1}^\m {N\choose
{k}}card(\adm(\g,k,k))\\
&&\hspace{-7.3cm}= \lim_{N\to\infty}\frac{1}{N^{\m}} \sum_{\g\in
PP(m)}V(\g,b_1,\dots,b_{m-1}) {N\choose
{\m}}card(\adm(\g,\m,\m))\\
&&\hspace{-7.3cm}=\frac{1}{\left(\m\right)!}\sum_{\g\in
PP(m)}V(\g,b_1,\dots,b_{m-1})card(\adm(\g,\m,\m))\\
\end{eqnarray*}

since

\[\lim_{N\to\infty}\frac{1}{N^\m}{N\choose{k}}=\left\{ \begin{array}{clcr}
0& \text{if} & k<\m\\
& & \\ \frac{1}{\left(\m\right)!}& \text{if} &
k=\m\end{array}\right.\]

For the last part, note that $V(\g,b_1,\dots,b_{m-1})$ is computed
iterating the result from Lemma \ref{5lemma1}, so for $\g\in PP(m)$
it depends only on the moments of order $2$ of $X_i$ ($i\geq1$).
\end{proof}

In the following paragraph we will suppose, without loss of
generality, that $\gB$ is unital.
\begin{cor}
The functional $\nu$ is positive if and only
 $\Phi(X_kb^*bX_k)\geq0$ for any $b\in\gB$ ($k\geq1$).
\end{cor}

\begin{proof}
One implication is trivial: if $\nu\geq0$, then
$\Phi(X_kb^*bX_k)=\nu\left((b\xi)^*b\xi\right)\geq0$.

 For the other implication, consider the set of symbols
 $\{\zeta_i\}_{i\geq1}$ and the linear spaces

 \[\gB \zeta_i\gB=\{b_1\zeta_ib_2,\ b_1,b_2\in\gB\}\]

\noindent with the $\gB$-bimodule structure given by

\[a_1(b_1\zeta_ib_2)a_2=(a_1b_1)\zeta_i(b_2a_2)\ \text{for any}\ a_1,a_2,b_1,b_2\in\gB\]

 \noindent and with the $\gB$-valued pairing $\langle,\rangle$ given by
 \begin{equation*}
 \langle a\zeta_i,b\zeta_i\rangle=
 \nu(\xi a^*b\xi)
\end{equation*}
\noindent for any $a,b\in\gB$.

The pairing $\langle, \rangle$ is positive, since $\nu(\xi
b^*b\xi)\geq0$ for any $b\in\gB$.

 Let $\cE$ be the weakly monotone product of
 $\{\gB\zeta_i\gB\}_{i\geq1}$. As shown in part \ref{weaklymonotone},
  the mappings $G(\zeta_i)$ form a
 monotonic independent family in $\mathcal{L}(\cE)$,
 therefore, from \ref{CLT}, one has that
 \[\nu(p^*(\xi)p(\xi))
 =\lim_{N\to \infty}\langle 1,p(\frac{G(\zeta_1)+\dots G(\zeta_N)}{\sqrt{N}})^*
 p(\frac{G(\zeta_1)+\dots G(\zeta_N)}{\sqrt{N}})1\rangle\]

\noindent for any $p(\xi)\in\gB\langle\xi\rangle$.

The conclusion follows from the positivity of the functional
$\langle1,\cdot1\rangle$.
\end{proof}

 \section{Positivity results and connection to operator-valued conditionally free products}
\begin{defn}
 \emph{Let $\gA_1, \gA_2$ be two algebras containing the subalgebra
 $\gB$ such that $\gA_1$ has the decomposition $\gA_1=\gB\oplus\gA_1^0$
for $\gA_1^0$ a subalgebra of $\gA$ which is also a $\gB$-algebra.
If $\Phi_1$, $\Phi_2$ are conditional expectations,
$\Phi_j:\gA_j\lra\gB, j=1,2$, we define $\Phi = \Phi_1\rhd\Phi_2$,
the monotonic product of $\Phi_1$ and
 $\Phi_2$ to be the conditional expectation on the algebraic free
 product with amalgamation over $\gB$,}
 \[\gA=\gA_1\ast_{\gB}\gA_2\]
 \emph{given by}
 \begin{eqnarray*}
 \Phi(\alpha a_1ba_2\beta)&=&\Phi(\alpha a_1\Phi_2(b)a_2\beta)\\
\Phi(ba_2\beta)&=&\Phi_2(b)\Phi(a_2\beta)\\
\Phi(\alpha a_1b)&=&\Phi(\alpha a_1)\Phi_2(b)
 \end{eqnarray*}
 \emph{for all}
 $a_1,a_2\in\gA_1^0, b\in\gA_2$ \emph{and} $\alpha, \beta\in\gA$.
\end{defn}
The map $\Phi$ is well-defined, since any element of $\gA$ can be
written as a sum of finite products in which the elements from
$\gA_1^0$ and $\gA_2$ and the conditions above imply
\begin{equation*}
\Phi(b_oa_1b_1\dots a_nb_n)=\Phi_1(\Phi_2(b_0)a_1\Phi_2(b_1)\dots
a_n\Phi_2(b_n))
\end{equation*}
for all $b_0,b_1,\dots,b_n\in\gB_2, a_1,\dots, a_n\in\gA_1^0$, and
all the analogues for the other types of such products.
\begin{prop}\label{monotposit}
If, in the above setting, $\gA_1, \gA_2$ are $*$-algebras, $\gB$ is
a $C^*$-algebra,  and $\Phi_1, \Phi_2$ are positive (i.e.
$\Phi_j(a^*a)\geq0$, for all $a\in\gA_j, j=1,2$), then
$\Phi_1\rhd\Phi_2$ is also positive.
\end{prop}
\begin{proof}
First, remember that the positivity of the conditional expectations
$\Phi_j$ implies that $\Phi_j(x^*)=\left(\Phi_j(x)\right)^*$, for
all $x\in\gA_j$.

Also, the map $\Phi_2$ is completely positive, and for any
$b_1,\dots b_n\in\gA_2$, the element
$(\Phi_2(b_i^*b_j))_{i,j=1}^n\in M_n(\gB)$ is positive (see
\cite{speicher1}, Section 3.5).

We have to show that
\[\Phi(a^*a)\geq0\ \ \ \ \ \
 \text{for all}\ a\in\gA_1\ast_{\gB}\gA_2.\]
 Any such $a$ can be written as a finite sum of elements of the form
 $b_0a_1b_1\dots a_nb_n$ with
 $b_0,\dots,b_n\in\gA_2,\ a_1,\dots,a_n\in\gA_1^0$.
 Hence:
 \begin{eqnarray*}
 \Phi(a^*a)&=&
 \Phi\left((\sum_{i=1}^Nb_{0,i}a_{1,i}b_{1,i}\dots
 a_{n(i),i}b_{n(i),i})^*(\sum_{i=1}^Nb_{0,i}a_{1,i}b_{1,i}\dots
 a_{n(i),i}b_{n(i),i})\right)\\
 &=&
 \Phi\left(\sum_{i,j=1}^Nb_{n(i),i}^*a_{n(i),i}^*\dots
 a_{1,i}^*b_{0,i}^*b_{0,j}a_{1,j}b_{1,j}\dots
 a_{n(j),j}b_{n(j),j}\right)\\
 &=&
 \sum_{i,j=1}^N \Phi_1\left(
\Phi_2(b_{n(i),i}^*)a_{n(i),i}^*\dots
 a_{1,i}^*\Phi_2(b_{0,i}^*b_{0,j})a_{1,j}\Phi_2(b_{1,j})\dots
 a_{n(j),j}\Phi_2(b_{n(j),j})
 \right)
 \end{eqnarray*}

 Since $\left(\Phi_2(b_{0,i}^*b_{0,j})\right)_{i,j=1}^N\in M_N(\gB)\subset M_N(\gA_1)$ is positive,
 there exists a matrix $T\in M_N(\gA_1)$ such that
 \[\left(\Phi_2(b_{0,i}^*b_{0,j})\right)_{i,j=1}^N=T^*T.\]

 With the notation
 \[a_i=a_{1,j}\Phi_2(b_{1,j})\dots
 a_{n(j),j}\Phi_2(b_{n(j),j}\in\gA_1\]
 \noindent we obtain:
 \begin{eqnarray*}
  \Phi(a^*a)
  &=&
  \Phi_1\left([(a_1\dots a_N)^*]T^*T(a_1\dots a_N)\right)\\
  &\geq&0
 \end{eqnarray*}
\end{proof}

Let $\gB\langle\xi,\xi^*\rangle$ be the $\ast$-algebra of
polynomials in $\xi$ and $\xi^*$ described in Section 2. For
$X\in\gA$, consider $\gA_X$ the $\ast$-subalgebra of $\gA$ generated
by $X$ and $\gB$. Define the mapping $\tau_X:\gB\langle\xi,
\xi^*\rangle\lra\gA_X$ to be the algebra $\ast$-homomorphism given
by $\tau_X(\xi)=X$ and the $\gB$-functional $\nu_X:\gB\langle\xi,
\xi^*\rangle\lra\gB$ to be given by $\nu_X=\Phi\circ\tau_X$.
\begin{cor}
If $X,Y\in\gA$ are monotonically independent over $\gB$ and $\nu_X,
\nu_Y$ are positive, then $\nu_{Z}$ is also positive for any element
$Z$ in the $\ast$-algebra generated by $X$ and $Y$. In particular
$\nu_{X+Y}$ and $\nu_{XY}$ are positive.
\end{cor}
\begin{proof}
 Consider $Z=Z(X,Y)$ a polynomial in $X$ and $Y$.
Since the maps
\begin{eqnarray*}
&&\nu_X:\gB\langle\xi_1,\xi_1^*\rangle\lra\gB\hspace{4cm}\\
&&\nu_Y:\gB\langle\xi_2, \xi_2^*\rangle\lra\gB
\end{eqnarray*}
are positive, from Proposition \ref{monotposit} so is

\[\nu_x\rhd\nu_Y:\gB\langle\xi_1,\xi_1^*\rangle\ast_\gB
\gB\langle\xi_2, \xi_2^*\rangle
=\gB\langle\xi_1,\xi_1^*,\xi_2,\xi_2^*\rangle\lra\gB\]

Remark also that
\[i_Z:\gB\langle\xi,\xi^*\rangle\ni f(\xi)\mapsto f(Z(\xi_1,\xi_2))
\in\gB\langle\xi_1,\xi_1^*\rangle\ast_\gB\gB\langle\xi_2,\xi_2^*\rangle\]
is a positive $\gB$-functional.

The conclusion follows from the fact that the monotonic independence
over $\gB$ of $X$ and $Y$ is equivalent to
\[\nu_{Z}=(\nu_X\rhd\nu_Y)\circ i_Z.\]

\end{proof}

 \begin{lemma}\label{mixprod}
 Let $\gA_1$, $\gA_2$  be two $\ast$-algebras containing the
 $C^*$-algebra $\gB$.
  and
 $\Phi_j:\gA_j\lra\gB, j=1,2$
 positive conditional
 expectations. Let $a_1,\dots,a_n\in\gA_1$, $a_{n+1},\dots,a_{n+m}\in\gA_2$
 and\\
  $A=(A_{i,j})\in M_{n+m}(\gB)$ be the matrix with the entries
 \[A_{i,j}=\left\{
 \begin{matrix}
 \Phi_1(a_i^*a_j)& \text{if}& i,j\leq n
  \hspace{1.2cm} \\
 \Phi_1(a_i^*)\Phi_2(a_j)& \text{if}& i\leq n, j>n\\
\Phi_2(a_i^*)\Phi_1(a_j)& \text{if}& i>n, j\leq n\\
\Phi_2(a_i^*a_j)& \text{if}& i,j> n\hspace{1.2cm}\\
 \end{matrix}\right.\]

 Then $A$ is positive.
 \end{lemma}
 \begin{proof}
 As shown in \cite{speicher1}, Theorem 3.5.6,  the
 $\gB$-functional $\Phi_1\ast_{\gB}\Phi_2$ is completely
 positive
 on $\gA_1\ast_\gB\gA_2$. Also, note that
 $A_{i,j}=\left((\Phi_1\ast_{\gB}\Phi_2)(a_i^*a_j)\right)$ for all $1\leq i,j=1\leq{n+m}$,
 and the conclusion follows from \cite{speicher1}, Lemma
 3.5.2.
  \end{proof}

 Consider now $\gA_1$, $\gA_2$  two $\ast$-algebras over the $C^*$-algebra
 $\gB$,
 each endowed with two positive conditional
 expectations
 $\Phi_j, \Psi_j:\gA_j\lra\gB$, $j=1,2$.
We define $(\gA, \Phi, \Psi)$, the \emph{conditionally free product
with amalgamation over $\gB$} of the triples $(\gA_1, \Phi_1,
\Psi_1)$ and $(\gA_2, \Phi_2, \Psi_2)$, by:
\begin{enumerate}
\item[(1)]$\gA=\gA_1\ast_\gB\gA_2$
\item[(2)]$\Psi=\Psi_1\ast_\gB\Psi_2$
and
 $ \Phi=\Phi_1\ast_{(\Psi_1,\Psi_2)}\Phi_2$
, i.e. the functionals $\Psi$ and $\Phi$ are determined by the
relations
\begin{eqnarray*}
\Psi(a_1a_2\dots a_n)&=&0\\
\Phi(a_1a_2\dots a_n)&=&\Phi(a_1)\Phi(a_2)\dots \Phi(a_n)
\end{eqnarray*}
 for any $a_i\in\gA_{\varepsilon(i)},
\varepsilon(i)\in\{1,2\},$\ such that
$\varepsilon(1)\neq\varepsilon(2)\neq\dots\neq \varepsilon(n)$\ and
$\Psi_{\varepsilon(i)}(a_i)=0$.
\end{enumerate}
\begin{thm}
In the above setting, $\Phi$ and $\Psi$ are positive
$\gB$-functionals.
\end{thm}\label{condposit}
\begin{proof}

The positivity of $\Psi$ is proved in \cite{speicher1}, Theorem
3.5.6.

For the positivity of $\Phi$ we have to show that $\Phi(a^*a)\geq 0$
for any $a\in\gA$.

Since any element of $\gA$ can be written
\begin{eqnarray*}
a&=& \sum_{k=1}^{N}s_{1,k}\dots s_{n(k),k}, \\
&& \ \ \ \ \ \ \ \
\text{where}\ s_{j,k}\in\gA_{\varepsilon(j,k)}\
\varepsilon(1,k)\neq\varepsilon(2,k)\neq\dots\neq \varepsilon(n(k),k)\\
&=&\sum_{k=1}^{N}\prod_{j=1}^{n(k)}
(s_{(j,k)}-\Psi(s_{(j,k))}+\Psi(s_{(j,k))}),
\end{eqnarray*}
we can consider $a$ of the form
\[a=\alpha+\sum_{k=1}^{N}a_{1,k}\dots a_{n(k),k}\]
such that
\begin{eqnarray*}
\alpha&\in&\gB\\
a_{j,k}&\in&\gA_{\varepsilon(j,k)},\
\varepsilon(1,k)\neq\varepsilon(2,k)\neq\dots\neq
\varepsilon(n(k),k)\\
\Psi_{\varepsilon(j,k)}(a_{j,k})&=&0.
\end{eqnarray*}
Therefore
\begin{eqnarray*}
\Phi(a^*a)
&=&
\Phi\left( \left(\alpha+\sum_{k=1}^{N}a_{1,k}\dots
a_{n(k),k}\right)^*\left(\alpha+\sum_{k=1}^{N}a_{1,k}\dots
a_{n(k),k}\right)\right)\\
&=&
\Phi\left(\alpha^*\alpha+\alpha^*\left(\sum_{k=1}^{N}a_{1,k}\dots
a_{n(k),k}\right)+\left(\sum_{k=1}^{N}a_{1,k}\dots
a_{n(k),k}\right)^*\alpha+\right.\\
&&\left.\left(\sum_{k=1}^{N}a_{1,k}\dots
a_{n(k),k}\right)^*\left(\sum_{k=1}^{N}a_{1,k}\dots
a_{n(k),k}\right)\right)\\
&=& \Phi(\alpha^*\alpha)+\sum_{k=1}^{N}\Phi(\alpha^* a_{1,k}\dots
a_{n(k),k})+\sum_{k=1}^{N}\Phi( a_{n(k),k}^*\dots a_{1,k}^*\alpha)\\
&&+ \sum_{k,l=1}^{N}\Phi( a_{n(k),k}^*\dots a_{1,k}^* a_{1,l}\dots
a_{n(l),l})\\
\end{eqnarray*}

Using the definition of the conditionally free product with
amalgamation over $\gB$ and that
$\Psi_{\varepsilon(j,k)}(a_{j,k})=0$ for all $j,k$, the above
relation becomes:

\begin{eqnarray*}
\Phi(a^*a) &=&
 \Phi(\alpha^*\alpha)+\sum_{k=1}^{N}\Phi(\alpha^*
 a_{1,k})\Phi(a_{2,k})\dots\Phi(a_{n(k),k})\\
 &&+\sum_{k=1}^{N}\Phi( a_{n(k),k})^*\dots \Phi(a_{2,k}^*)\Phi(a_{1,k}^*\alpha)\\
 &&+\sum_{k,l=1}^{N}\left(\Phi( a_{n(k),k})^*\dots\Phi(a_{2,k}^*)\right)\Phi(a_{1,k}^*
 a_{1,l})\Phi(a_{2,l})\dots\Phi(a_{n(l),l})\\
 \end{eqnarray*}
that is
\begin{eqnarray*}
\Phi(a^*a)
 &=&
 \Phi(\alpha^*\alpha)+\sum_{k=1}^{N}\Phi(\alpha^*
 a_{1,k})\left[\Phi(a_{2,k})\dots\Phi(a_{n(k),k})\right]\\
 &&+\sum_{k=1}^{N}\left[\Phi(a_{2,k})\dots \Phi( a_{n(k),k})\right]^*\Phi(a_{1,k}^*\alpha)\\
 &&+\sum_{k,l=1}^{N}\left[\Phi(a_{2,k})\dots \Phi( a_{n(k),k})\right]^*\Phi(a_{1,k}^*
 a_{1,l})\left[\Phi(a_{2,l})\dots\Phi(a_{n(l),l})\right]\\
\end{eqnarray*}

Denote now $a_{1,N+1}=\alpha$ and
$\beta_k=\Phi(a_{2,k})\dots\Phi(a_{n(k),k})$.

 From Lemma \ref{mixprod},
 the matrix $S=\left(\Phi(a_{1,i}^*a_{1,j})_{i,j=1}^{N+1}\right)$
 is positive in $M_{N+1}(\gB)$, therefore
 \[S=T^*T,\  \text{for some}\ T\in M_{N+1}(\gB)\]

 The identity for $\Phi(a^*a)$ becomes:
 \begin{eqnarray*}
\Phi(a^*a)&=&(\beta_1,\dots,\beta_N,1)^*T^*T(\beta_1,\dots,\beta_N,1)\\
&\geq&0,
 \end{eqnarray*}

as claimed.

\end{proof}

Suppose now that the $\ast$-algebra $\gA_1$ has the decomposition
$\gA_1=\gB\oplus \gA_1^{0}$, such that $\gA_1^{0}$ is a
$\ast$-subalgebra of $\gA_1$ which is also a $\gB$-algebra. Define
the $\gB$-valued conditional expectation
\[\delta:\gA_1\ni(\lambda+a_0)\mapsto\lambda\in\gB\]
for all $a_0\in\gA_1^0$.

\begin{thm}
With the notations above,
\[\Phi_1\triangleright\Phi_2=\Phi_1\ast_{(\delta,\Phi_2)}\Phi_2\]
\end{thm}
\begin{proof}
First remark that $\delta(a)=0$ implies $a\in\gA_1^0$, from the
definition of $\delta$.

 For
$\varepsilon(1),\dots,\varepsilon(n)\in\{1,2\}$
such that
$\varepsilon(1)\neq\dots\neq\varepsilon(n)$
 and
$a_j\in\gA_{\varepsilon(j)}$ such that $\delta(a_j)=0$ if
$\varepsilon(j)=1$ and $\Phi_2(a_j)=0$ if $\varepsilon(j)=2$, one
has
($\chi_{\gA_j}$ denotes the characteristic function of $\gA_j$):
\begin{eqnarray*} (\Phi_1\triangleright\Phi_2)(a_1\dots a_n) &=&
\Phi_1\left(\prod_{j=1}^n\left[
\chi_{\gA_1}(a_j)+\Phi_2(\chi_{\gA_2}(a_j))\right]\right)\\
&=&\Phi_1\left(\prod_{j=1}^n
\chi_{\gA_1}(a_j)\right)\\
 &=&0\ \ \text{since there is at least one}\ a_j\in\gA_2
\end{eqnarray*}

The conclusion follows from the above equality, since the
conditional expectation
 $\Phi_1\ast_{\delta,\Phi_2}\Phi_2$
 is generated by
\[ (\Phi_1\ast_{(\delta,\Phi_2)}\Phi_2)(a_1\dots a_n)=0.\]
\end{proof}

\textbf{Acknowledgements.} The true motivation for this paper was a
discussion during my stay in Wroclaw (with the generous support of
EU Network QP-Applications contract HPRN-CT-2002-00729) with Janusz
Wysocza\`{n}ski, who brought to my attention the reference
\cite{janusz}, for which I am very thankful. I am also indebted to
Wojtech M{\l}otkowski for the reference \cite{mlotk1}, but
especially to Hari Bercovici for his priceless help in conceiving
and writing this material.

\end{document}